\documentclass[showpacs,aps,prb]{revtex4}

\usepackage{amsmath,amssymb,amsthm}
\usepackage{mathtools}
\usepackage{xspace}
\usepackage[usenames,dvipsnames]{xcolor}
\usepackage[colorlinks,hyperindex,breaklinks]{hyperref}
\usepackage{braket}
\usepackage{hyphenat}
\usepackage{bbold}
\usepackage{enumerate}
\usepackage{tabularx}
\usepackage[hang,flushmargin]{footmisc}
\usepackage[british]{babel}

\usepackage{srcltx}
\usepackage{color}

\usepackage{graphicx}
\usepackage{amsmath}
\usepackage{amssymb}
\usepackage{bm}

\begin{document}

\title{On four families of power series involving harmonic numbers and central binomial coefficients}
\author{J. Braun}
\affiliation{Ludwig Maximilians-Universit{\"a}t, M{\"u}nchen, Germany}
\author{D. Romberger}
\affiliation{Fakultät IV, Abt. BWL, Hochschule Hannover}
\author{H. J. Bentz}
\affiliation{Institut f\"{u}r Mathematik und Informatik, Universit\"{a}t Hildesheim, Hildesheim, Germany}
\date{\today}

\begin{abstract}
We present several sequences involving harmonic numbers and the central binomial coefficients. The calculational
technique is consists of a special summation method that allows, based on proper two-valued integer functions, to
calculate different families of power series which involve odd harmonic numbers and central binomial coefficients.
Furthermore it is shown that based on these series a new type of nonlinear Euler sums that involve odd harmonic
numbers can be calculated in terms of zeta functions.
\end{abstract}
\maketitle

\section{Introduction}
The central binomial coefficients are closely related to the so called Catalan numbers. Many facts about these
coefficients and Catalan numbers can be found, for example, by Koshy \cite{koh08}. A variety of identities involving
central binomial coefficients has been collected by Gould \cite{gol79}. Also Riordan presents a large list of
references on Catalan numbers in his book \cite{rio79}. Further identities involving central Binomial coefficients
can be found in \cite{gen11}. We focus here on families of power series which involve odd harmonic numbers and
central binomial coefficients. Earlier papers which focused on similar power series were published by Zucker, Lehmer,
Boyadzhiev and others \cite{leh85,wei04,zuc85,han75,boy12,che16}. In this work we introduce four families of power
series involving inverse powers in combination with central binomial numbers as well as series involving inverse powers
multiplied by odd harmonic numbers in combination with central binomial coefficients.

\section{A first kind of inverse power series with central Binomial coefficients}

The first family of series is defined by the following equation:

\begin{eqnarray}
s(n) = \sum^{\infty}_{k=1}\frac{1}{k^n}
\left(\begin{array}{c}
                               2k \\
                               k
\end{array}\right)\frac{1}{4^k}~,
\end{eqnarray}
with $n \in \mathbb{N}$.

As an example, for n=1 and 2 the series are known from literature \cite{boy12}. It follows: 
\begin{eqnarray}
\sum^{\infty}_{k=1}\frac{1}{k}
\left(\begin{array}{c}
                               2k \\
                               k
\end{array}\right)\frac{1}{4^k} = 2ln(2)
\end{eqnarray}

\begin{eqnarray}
\sum^{\infty}_{k=1}\frac{1}{k^2}
\left(\begin{array}{c}
                               2k \\
                               k
\end{array}\right)\frac{1}{4^k} = \zeta(2) - 2\left( ln(2) \right)^2~.
\end{eqnarray}

To be able to compute the corresponding series for higher n values we need two identities in form of proper valued integer series.

\subsection{Lemma 1}

The following identity holds:
\begin{eqnarray}
f(k) = \sum^{\infty}_{i=1}\frac{1}{i+k}
\left(\begin{array}{c}
                               2i \\
                               i
               \end{array}\right) \frac{1}{4^i} =
\frac{1}{k}\left(\begin{array}{c}
                               2k \\
                               k
\end{array}\right)^{-1} 4^k - \frac{1}{k}
\end{eqnarray}

\subsection{Proof of Lemma 1}
For k=0 the result can be obtained from the following generating function \cite{boy12}: 
\begin{eqnarray}
\sum^{\infty}_{i=1}\frac{1}{i}
\left(\begin{array}{c}
                               2i \\
                               i
\end{array}\right) x^i = 2ln \left( \frac{2}{1+\sqrt{1-4x}} \right)~.
\end{eqnarray}
With x = $\frac{1}{4}$ one has f(0) = 2ln(2). For the case k=1 we start with:

\begin{eqnarray}
f(1) = \sum^{\infty}_{i=1}\frac{1}{i+1}
\left(\begin{array}{c}
                               2i \\
                               i
\end{array}\right) \frac{1}{4^i} \nonumber
\end{eqnarray}
\begin{eqnarray}
= \sum^{\infty}_{i=2}\frac{1}{i}
\left(\begin{array}{c}
                               2i-2 \\
                               i-1
\end{array}\right) \frac{1}{4^{i-1}} \nonumber
\end{eqnarray}
\begin{eqnarray}
= 2\sum^{\infty}_{i=1}\frac{1}{2i-1}
\left(\begin{array}{c}
                               2i \\
                               i
\end{array}\right) \frac{1}{4^{i}} - 1
\end{eqnarray}
This way it remains to compute the last expression. We define:
\begin{eqnarray}
u(i) =\sum^{\infty}_{i=1}\frac{h(i)}{2i-1}
\left(\begin{array}{c}
                               2i \\
                               i
\end{array}\right) \frac{1}{4^{i}}
\end{eqnarray}
with h(i) an arbitrary integer function. It follows first: 
\begin{eqnarray}
u(i) =2\sum^{\infty}_{i=1}\frac{h(i)}{i}
\left(\begin{array}{c}
                               2i-2 \\
                               i-1
\end{array}\right) \frac{1}{4^{i}} \nonumber
\end{eqnarray} 
\begin{eqnarray}
= \frac{1}{2}\sum^{\infty}_{i=0}\frac{h(i+1)}{i+1}
\left(\begin{array}{c}
                               2i \\
                               i
\end{array}\right) \frac{1}{4^{i}}
\end{eqnarray}
With the identity
\begin{eqnarray}
\frac{1}{i+1}
\left(\begin{array}{c}
                               2i \\
                               i
\end{array}\right) = 2
\left(\begin{array}{c}
                               2i \\
                               i
\end{array}\right) -
\left(\begin{array}{c}
                               2i+1 \\
                               i+1
\end{array}\right)
\end{eqnarray}
we get
\begin{eqnarray}
u(i) = \sum^{\infty}_{i=0}h(i+1)
\left(\begin{array}{c}
                               2i \\
                               i
\end{array}\right) \frac{1}{4^{i}} - \frac{1}{2}
\sum^{\infty}_{i=0}h(i+1)\left(\begin{array}{c}
                               2i+1 \\
                               i+1
\end{array}\right) \frac{1}{4^{i}} \nonumber
\end{eqnarray} 
\begin{eqnarray}
= \sum^{\infty}_{i=0}h(i+1)
\left(\begin{array}{c}
                               2i \\
                               i
\end{array}\right) \frac{1}{4^{i}} -2 
\sum^{\infty}_{i=1}h(i)\left(\begin{array}{c}
                               2i-1 \\
                               i
\end{array}\right) \frac{1}{4^{i}} \nonumber
\end{eqnarray} 
\begin{eqnarray}
= h(1) + \sum^{\infty}_{i=1}\left(h(i+1) - h(i)\right)
\left(\begin{array}{c}
                               2i \\
                               i
\end{array}\right) \frac{1}{4^{i}}
\end{eqnarray}.
With h(i) = 1 it follows u(1) = 1 and with this f(1) = 1. Now we calculate 
\begin{eqnarray}
f(2) = \sum^{\infty}_{i=1}\frac{1}{i+2}
\left(\begin{array}{c}
                               2i \\
                               i
\end{array}\right) \frac{1}{4^i} \nonumber
\end{eqnarray}
\begin{eqnarray}
= 4\sum^{\infty}_{i=2}\frac{1}{i+1}
\left(\begin{array}{c}
                               2i-2 \\
                               i-1
\end{array}\right) \frac{1}{4^i} \nonumber
\end{eqnarray}
\begin{eqnarray}
= \frac{2}{3}\sum^{\infty}_{i=1}\frac{1}{2i-1}
\left(\begin{array}{c}
                               2i \\
                               i
\end{array}\right) \frac{1}{4^i} + 
\frac{2}{3}\sum^{\infty}_{i=1}\frac{1}{i+1}
\left(\begin{array}{c}
                               2i \\
                               i
\end{array}\right) \frac{1}{4^i} - \frac{1}{2} 
\end{eqnarray}
with the result f(2) = $\frac{5}{6}$.
Repeating this procedure k-times we end up with the following equation:
\begin{eqnarray}
f(k) - \frac{2k-2}{2k-1} f(k-1) = \frac{1}{k(2k-1)}~.
\end{eqnarray}
This is a inhomogeneous difference equation of first order with a non-constant coefficient.
The solution of the homogenous equation is:
\begin{eqnarray}
f(k) = \frac{1}{2k}
\left(\begin{array}{c}
                               2k \\
                               k
\end{array}\right)^{-1}4^k
\end{eqnarray}
and a special solution of the inhomogeneous equation results to:
\begin{eqnarray}
f(k) = \frac{1}{2k}
\left(\begin{array}{c}
                               2k \\
                               k
\end{array}\right)^{-1}4^k
\sum^{k}_{i=1} \frac{2i}{i(2i-1)}
\left(\begin{array}{c}
                               2k \\
                               k
\end{array}\right) \frac{1}{4^k} \nonumber
\end{eqnarray}
\begin{eqnarray}
= \frac{1}{k}
\left(\begin{array}{c}
                               2k \\
                               k
\end{array}\right)^{-1}4^k
\left( 1 - 
\left(\begin{array}{c}
                               2k \\
                               k
\end{array}\right) \frac{1}{4^k} \right)~.  
\end{eqnarray}
Thus Lemma 1 is proved.

\subsection{Lemma 2}
It holds 
\begin{eqnarray}
f(k) = \sum^{\infty}_{i=1} \frac{1}{(i+k)^2}
\left(\begin{array}{c}
                               2i \\
                               i
\end{array}\right) \frac{1}{4^i} =
\left( \frac{1}{k^2} - 2ln(2)\frac{1}{k} + \frac{2h_k}{k} - \frac{H_k}{k} \right)
\left(\begin{array}{c}
                               2k \\
                               k
\end{array}\right)^{-1} 4^k
-\frac{1}{k^2}~, 
\end{eqnarray}
with $k \in \mathbb{N}$.

\subsection{Proof of Lemma 2}
For k=0 the result can be obtained from Boyadzhiev \cite{boy12}. It follows:
\begin{eqnarray}
\sum^{\infty}_{i=1}\frac{1}{i^2}
\left(\begin{array}{c}
                               2i \\
                               i
\end{array}\right) \frac{1}{4^k} = \zeta(2) - 2\left( ln(2)\right)^2~.
\end{eqnarray}
For the case k=1 we start with:
\begin{eqnarray}
f(1) = \sum^{\infty}_{i=1}\frac{1}{(i+1)^2}
\left(\begin{array}{c}
                               2i \\
                               i
\end{array}\right) \frac{1}{4^k} \nonumber
\end{eqnarray}
\begin{eqnarray}
= 4\sum^{\infty}_{i=2}\frac{1}{i^2}
\left(\begin{array}{c}
                               2i-2 \\
                               i-1
\end{array}\right) \frac{1}{4^k} \nonumber 
\end{eqnarray}
\begin{eqnarray}
= 2\sum^{\infty}_{i=2}\frac{1}{i(2i-1)}
\left(\begin{array}{c}
                               2i \\
                               i
\end{array}\right) \frac{1}{4^k} \nonumber
\end{eqnarray}
\begin{eqnarray}
= 4\sum^{\infty}_{i=2}\frac{1}{2i-1}
\left(\begin{array}{c}
                               2i \\
                               i
\end{array}\right) \frac{1}{4^k} - 
2\sum^{\infty}_{i=2}\frac{1}{i}
\left(\begin{array}{c}
                               2i \\
                               i
\end{array}\right) \frac{1}{4^k} - 1 = 3 - 4ln(2)~. 
\end{eqnarray}
Repeating this procedure k-times we end up with the following inhomogeneous difference equation:
\begin{eqnarray}
f(k) - \frac{2k-2}{2k-1} f(k-1) = \frac{4}{(2k-1)^2} - \frac{1}{k^2} - \frac{2}{(2k-1)^2}
\sum^{\infty}_{k=1}\frac{1}{k+n-1}
\left(\begin{array}{c}
                               2k \\
                               k
\end{array}\right)\frac{1}{4^k}~,
\end{eqnarray}
where the solution of the homogeneous equation is known from the proof of lemma 1. Thus it follows:
\begin{eqnarray}
f(k) = \frac{1}{2k}
\left(\begin{array}{c}
                               2k \\
                               k
\end{array}\right)^{-1}4^k
	\sum^{k}_{i=1} \left( \frac{4}{(2i-1)^2} - \frac{1}{i^2} - \frac{2}{(2i-1)^2}
\sum^{\infty}_{n=1}\frac{1}{n+i-1}
\left(\begin{array}{c}
                               2n \\
                               n
\end{array}\right)\frac{1}{4^n} \right) 2i
\left(\begin{array}{c}
                               2i \\
                               i
\end{array}\right)\frac{1}{4^i}~.
\end{eqnarray}
For k=1 the sum results to $f(1)=3-4ln(2)$. It is advantageous to start the summation with the index k=2 and to add the first term in
its explicit form. It follows then: 
\begin{eqnarray}
f(k) = \frac{1}{2k}
\left(\begin{array}{c}
                               2k \\
                               k
\end{array}\right)^{-1}4^k 
\left[3-4ln(2) + 
\sum^{k}_{i=2} \left( \frac{4}{(2i-1)^2} - \frac{1}{i^2} - \frac{2}{(2i-1)^2}
\sum^{\infty}_{n=1}\frac{1}{n+i-1}
\left(\begin{array}{c}
                               2n \\
                               n
\end{array}\right)\frac{1}{4^n} \right)
\left(\begin{array}{c}
                               2i \\
                               i
\end{array}\right)\frac{2i}{4^i} \right]  \nonumber
\end{eqnarray}
\begin{eqnarray}
 = \frac{1}{2k}
\left(\begin{array}{c}
                               2k \\
                               k
\end{array}\right)^{-1}4^k 
\left[3-4ln(2) + 
\sum^{k-1}_{i=1} \left( \frac{4}{(2i+1)^2} - \frac{1}{(i+1)^2} - \frac{2}{(2i+1)^2}
\sum^{\infty}_{n=1}\frac{1}{n+i}
\left(\begin{array}{c}
                               2n \\
                               n
\end{array}\right)\frac{1}{4^n} \right)
\left(\begin{array}{c}
                               2i \\
                               i
\end{array}\right)\frac{2i+1}{4^i} \right] \nonumber
\end{eqnarray}
\begin{eqnarray}
 = \frac{1}{2k}
\left(\begin{array}{c}
                               2k \\
                               k
\end{array}\right)^{-1}4^k 
\left[3-4ln(2) + 
\sum^{k-1}_{i=1} \left( \frac{4}{(2i+1)^2} - \frac{1}{(i+1)^2} - \frac{2}{i(2i+1)^2}
\left( \left(\begin{array}{c}
                               2i \\
                               i
\end{array}\right)^{-1}{4^i} - 1 \right) \right)
\left(\begin{array}{c}
                               2i \\
                               i
\end{array}\right)\frac{2i+1}{4^i} \right] \nonumber
\end{eqnarray}
\begin{eqnarray}
 = \frac{1}{2k}
\left(\begin{array}{c}
                               2k \\
                               k
\end{array}\right)^{-1}4^k 
\left[ 3-4ln(2) + 
\sum^{k-1}_{i=1} \left( \frac{4}{(2i+1)} - \frac{2i+1}{(i+1)^2} + \frac{2}{i(2i+1)} \right)
\left(\begin{array}{c}
                               2i \\
                               i
\end{array}\right)\frac{1}{4^i} - 2\sum^{k-1}_{i=1} \frac{1}{i(2i+1)} \right] \nonumber
\end{eqnarray}
\begin{eqnarray}
 = \frac{1}{2k}
\left(\begin{array}{c}
                               2k \\
                               k
\end{array}\right)^{-1}4^k 
\left[ 3-4ln(2) + 
\sum^{k-1}_{i=1} \frac{2}{i}
\left(\begin{array}{c}
                               2i \\
                               i
\end{array}\right)\frac{1}{4^i}-
\frac{1}{2} \sum^{k-1}_{i=1} \frac{1}{i+1}
\left(\begin{array}{c}
                               2i+2 \\
                               i+1
\end{array}\right)\frac{1}{4^i} - 2\sum^{k-1}_{i=1} \frac{1}{i(2i+1)} \right] \nonumber 
\end{eqnarray}
\begin{eqnarray}
 = \frac{1}{2k}
\left(\begin{array}{c}
                               2k \\
                               k
\end{array}\right)^{-1}4^k
\left[ 3-4ln(2) + 1 - \frac{2}{k}
\left(\begin{array}{c}
                               2k \\
                               k
\end{array}\right)\frac{1}{4^k}-
2\sum^{k-1}_{i=1} \frac{1}{i(2i+1)} \right] \nonumber
\end{eqnarray}
\begin{eqnarray}
 = \frac{1}{2k}
\left(\begin{array}{c}
                               2k \\
                               k
\end{array}\right)^{-1}4^k
	\left( 3-4ln(2) + 1 - 2H_k + 4 h_k + \frac{2}{k} - 4 \right) -\frac{1}{k^2} \nonumber
\end{eqnarray}
\begin{eqnarray}
=\left( \frac{1}{k^2} - 2ln(2)\frac{1}{k} + \frac{2h_k}{k} - \frac{H_k}{k} \right)
\left(\begin{array}{c}
                               2k \\
                               k
\end{array}\right)^{-1} 4^k
-\frac{1}{k^2}~. 
\end{eqnarray}
Thus lemma 2 is proved.

Having finished the proof for lemma 2 we can start to compute the corresponding series for n-values higher than 2.
It follows first for n=3: 
\begin{eqnarray}
\frac{1}{i(k+i)^2} = \frac{1}{k^2}\frac{1}{i} - \frac{1}{k^2}\frac{1}{k+i} - \frac{1}{k}\frac{1}{(k+i)^2}~.
\end{eqnarray}
Now we write:
\begin{eqnarray}
\sum^{\infty}_{i=1}\frac{1}{i}
\left(\begin{array}{c}
                               2i \\
                               i
\end{array}\right)\frac{1}{4^i}
	\sum^{\infty}_{k=1}\frac{1}{(k+i)^2}
\left(\begin{array}{c}
                               2k \\
                               k
\end{array}\right)\frac{1}{4^k} \nonumber
\end{eqnarray}
\begin{eqnarray}
= \frac{1}{2}\sum^{\infty}_{k=1}\frac{1}{k^2}
\left(\begin{array}{c}
                               2k \\
                               k
\end{array}\right)\frac{1}{4^k}
\sum^{\infty}_{i=1}\frac{1}{i}
\left(\begin{array}{c}
                               2i \\
                               i
\end{array}\right)\frac{1}{4^i}
- \frac{1}{2}\sum^{\infty}_{k=1}\frac{1}{k^2}
\left(\begin{array}{c}
                               2k \\
                               k
\end{array}\right)\frac{1}{4^k}
\sum^{\infty}_{i=1}\frac{1}{k+i}
\left(\begin{array}{c}
                               2i \\
                               i
\end{array}\right)\frac{1}{4^i} \nonumber
\end{eqnarray}
\begin{eqnarray}
= \frac{1}{2}\sum^{\infty}_{k=1}\frac{1}{k^2}
\left(\begin{array}{c}
                               2k \\
                               k
\end{array}\right)\frac{1}{4^k}
\sum^{\infty}_{i=1}\frac{1}{i}
\left(\begin{array}{c}
                               2i \\
                               i
\end{array}\right)\frac{1}{4^i}
- \frac{1}{2}\sum^{\infty}_{k=1}\frac{1}{k^2}
\left(\begin{array}{c}
                               2k \\
                               k
\end{array}\right)\frac{1}{4^k}\frac{1}{k} \left(
\left(\begin{array}{c}
                               2k \\
                               k
\end{array}\right)^{-1}4^k -1  \right) \nonumber
\end{eqnarray}
\begin{eqnarray}
= ln(2)\zeta(2) - 2\left( ln(2) \right)^3 - \frac{1}{2}\zeta(3) +\frac{1}{2}
\sum^{\infty}_{k=1}\frac{1}{k^3}
\left(\begin{array}{c}
                               2k \\
                               k
\end{array}\right)\frac{1}{4^k} \nonumber
\end{eqnarray}
\begin{eqnarray}
= 2\sum^{\infty}_{k=1}\frac{h_k}{k^2} - \sum^{\infty}_{k=1}\frac{H_k}{k^2} + \zeta(3) -2ln(2)\zeta(2) -
\sum^{\infty}_{k=1}\frac{1}{k^3}
\left(\begin{array}{c}
                               2k \\
                               k
\end{array}\right)\frac{1}{4^k}~.
\end{eqnarray}
From literature \cite{zeh07,ade16} it follows: 
\begin{eqnarray}
\sum^{\infty}_{k=1}\frac{H_k}{k^2}=2\zeta(3)~,
\end{eqnarray}
and 
\begin{eqnarray}
\sum^{\infty}_{k=1}\frac{h_k}{k^2}=\frac{7}{4}\zeta(3)~,
\end{eqnarray}
and with this
\begin{eqnarray}
\sum^{\infty}_{k=1}\frac{1}{k^3}
\left(\begin{array}{c}
                               2k \\
                               k
\end{array}\right)\frac{1}{4^k} = 2\zeta(3) - 2ln(2)\zeta(2) + \frac{4}{3} \left( ln(2) \right)^3~.
\end{eqnarray}
In the case that n=4 we start with the following expression: 
\begin{eqnarray}
\frac{1}{i^2(k+i)^2} = \frac{1}{k^2}\frac{1}{i^2} - \frac{2}{k^3}\frac{1}{i} + \frac{2}{k^3}\frac{1}{k+i} 
+ \frac{1}{k^2}\frac{1}{(k+i)^2}~.
\end{eqnarray}
Now we write:
\begin{eqnarray}
\sum^{\infty}_{i=1}\frac{1}{i^2}
\left(\begin{array}{c}
                               2i \\
                               i
\end{array}\right)\frac{1}{4^i}
	\sum^{\infty}_{k=1}\frac{1}{(k+i)^2}
\left(\begin{array}{c}
                               2k \\
                               k
\end{array}\right)\frac{1}{4^k} - 
\sum^{\infty}_{k=1}\frac{1}{k^2}
\left(\begin{array}{c}
                               2k \\
                               k
\end{array}\right)\frac{1}{4^k}
	\sum^{\infty}_{i=1}\frac{1}{(k+i)^2}
\left(\begin{array}{c}
                               2i \\
                               i
\end{array}\right)\frac{1}{4^i} =  \nonumber 
\end{eqnarray}
	
\begin{eqnarray}
\sum^{\infty}_{k=1}\frac{1}{k^2}
\left(\begin{array}{c}
                               2k \\
                               k
\end{array}\right)\frac{1}{4^k}
\sum^{\infty}_{i=1}\frac{1}{i^2}
\left(\begin{array}{c}
                               2i \\
                               i
\end{array}\right)\frac{1}{4^i}
- \sum^{\infty}_{k=1}\frac{2}{k^3}
\left(\begin{array}{c}
                               2k \\
                               k
\end{array}\right)\frac{1}{4^k}
\sum^{\infty}_{i=1}\frac{1}{i}
\left(\begin{array}{c}
                               2i \\
                               i
\end{array}\right)\frac{1}{4^i} +
\sum^{\infty}_{k=1}\frac{2}{k^3}
\left(\begin{array}{c}
                               2k \\
                               k
\end{array}\right)\frac{1}{4^k}
\sum^{\infty}_{i=1}\frac{1}{i+k}
\left(\begin{array}{c}
                               2i \\
                               i
\end{array}\right)\frac{1}{4^i}~.
\end{eqnarray}
The left side of equation (27) is zero. So we can write:
\begin{eqnarray}
\sum^{\infty}_{k=1}\frac{2}{k^3}
\left(\begin{array}{c}
                               2k \\
                               k
\end{array}\right)\frac{1}{4^k}
\sum^{\infty}_{i=1}\frac{1}{i+k}
\left(\begin{array}{c}
                               2i \\
                               i
\end{array}\right)\frac{1}{4^i} \nonumber
\end{eqnarray}
\begin{eqnarray}
=\sum^{\infty}_{k=1}\frac{2}{k^4}
\left(\begin{array}{c}
                               2k \\
                               k
\end{array}\right)\frac{1}{4^k} \left(
\left(\begin{array}{c}
                               2k \\
                               k
\end{array}\right)^{-1}4^k - 1 \right) \nonumber
\end{eqnarray}
\begin{eqnarray}
= 4ln(2)\left( 2\zeta(3) - 2ln(2)\zeta(2) + \frac{4}{3} \left( ln(2) \right)^3 \right) 
- \left( \zeta(2) - 2\left(ln(2)\right)^2 \right)^2
\end{eqnarray}
Thus we get after elementary manipulations: 
\begin{eqnarray}
\sum^{\infty}_{k=1}\frac{1}{k^4}
\left(\begin{array}{c}
                               2k \\
                               k
\end{array}\right)\frac{1}{4^k} = \frac{9}{4}\zeta(4) - 4ln(2)\zeta(3) + 2\left( ln(2)\right)^2 \zeta(2) - 
\frac{2}{3} \left( ln(2) \right)^4
\end{eqnarray}
This procedure works analogously for higher n values, where for odd powers of n the following harmonic series are needed:
\begin{eqnarray}
\sum^{\infty}_{k=1}\frac{H_k}{k^{2n}}~,
\end{eqnarray}
and 
\begin{eqnarray}
\sum^{\infty}_{k=1}\frac{h_k}{k^{2n}}~.
\end{eqnarray}
These series are known from literature \cite{sit85}. As a consequence the corresponding type of inverse power series can be
calculated for all $n \in \mathbb{N}$ recursively in terms of zeta functions and natural logarithms. As an example the series
for n=5,6,7 are shown below.

\begin{widetext}
\begin{eqnarray}
\sum^{\infty}_{k=1}\frac{1}{k^5}
\left(\begin{array}{c}
                               2k \\
                               k
\end{array}\right)\frac{1}{4^k} = 6\zeta(5)-2\zeta(2)\zeta(3) - \frac{9}{2}ln(2)\zeta(4) + 4\left( ln(2)\right)^2 \zeta(3)
- \frac{4}{3} \left( ln(2) \right)^3 \zeta(2) + \frac{4}{15} \left( ln(2) \right)^5
\end{eqnarray}

\begin{eqnarray}
\sum^{\infty}_{k=1}\frac{1}{k^6}
\left(\begin{array}{c}
                               2k \\
                               k
\end{array}\right)\frac{1}{4^k} &=& \frac{79}{16}\zeta(6)-12ln(2)\zeta(5) + 4ln(2)\zeta(2)\zeta(3) - 2\left(\zeta(3) \right)^2 +
\frac{9}{2}\left(ln(2)\right)^2 \zeta(4) - \frac{8}{3}\left( ln(2)\right)^3 \zeta(3) \nonumber \\
&+& ~\frac{2}{3}\left( ln(2)\right)^4 \zeta(2) - \frac{4}{45} \left( ln(2) \right)^6
\end{eqnarray}

\begin{eqnarray}
\sum^{\infty}_{k=1}\frac{1}{k^7}
\left(\begin{array}{c}
                               2k \\
                               k
\end{array}\right)
\frac{1}{4^k} &=& 18\zeta(7)-\frac{79}{8}ln(2)\zeta(6)-6\zeta(2)\zeta(5)+12\left(ln(2)\right)^2\zeta(5)
-\frac{9}{2}\zeta(3)\zeta(4)-3\left(ln(2)\right)^3\zeta(4) \nonumber \\
&+& 4ln(2)\zeta(3)^2-4\left(ln(2)\right)^2\zeta(2)\zeta(3)+\frac{4}{3}\left(ln(2)\right)^4\zeta(3)
-\frac{4}{15}\left(ln(2)\right)^5\zeta(2)+\frac{8}{315}\left(ln(2)\right)^7
\end{eqnarray}
\end{widetext}
It remains here to mention that for all powers in $n \in \mathbb{N}$ the corresponding series can be
explicitly calculated by combinations of zeta functions and logarithmic functions.

\section{A Second kind of inverse power series with central Binomial coefficients}

The second family of inverse power series is defined as follows:
\begin{eqnarray}
l(n) = \sum^{\infty}_{k=1}\frac{1}{(2k+1)^{n}}
\left(\begin{array}{c}
                               2k \\
                               k
\end{array}\right)\frac{1}{4^k}~,
\end{eqnarray}
with $n \in \mathbb{N}$. For $n\le 3$ the series are known from literature \cite{zuc85}
\begin{eqnarray}
\sum^{\infty}_{i=1} \frac{1}{(2i+1)}
\left(\begin{array}{c}
                               2i \\
                               i
               \end{array}\right) \frac{1}{4^i} 
= \frac{\pi}{2} - 1
\end{eqnarray}

\begin{eqnarray}
\sum^{\infty}_{i=1} \frac{1}{(2i+1)^2}
\left(\begin{array}{c}
                               2i \\
                               i
               \end{array}\right) \frac{1}{4^i} 
= \frac{\pi}{2}~ln(2) - 1
\end{eqnarray}

\begin{eqnarray}
\sum^{\infty}_{i=1} \frac{1}{(2i+1)^3}
\left(\begin{array}{c}
                               2i \\
                               i
               \end{array}\right) \frac{1}{4^i} 
= \frac{\pi}{8} \zeta(2) + \frac{\pi}{4}~(ln(2))^2 - 1
\end{eqnarray}

To be able to compute the corresponding series for higher n values we need the following identity.

\section{Theorem 1}
It holds:
\begin{eqnarray}
\sum^{\infty}_{k=1}\frac{1}{(2k+1)^{n+1}}
\left(\begin{array}{c}
                               2k \\
                               k
\end{array}\right)\frac{1}{4^k} = \frac{(-)^n}{n!} \int_0^{\pi/2} \left( lnsin(x) \right)^n dx - 1
\end{eqnarray}

\subsection{Proof of Theorem 1}
We start with the Taylor series:
\begin{eqnarray}
\frac{1}{\sqrt{1-x^2}} = \sum^{\infty}_{k=0}
\left(\begin{array}{c}
                               2k \\
                               k
\end{array}\right)\frac{1}{4^k} x^{2k}	
\end{eqnarray}
and integrate this series in the interval [0,1]:
\begin{eqnarray}
\int^1_0 \frac{1}{\sqrt{1-x^2}} = arcsin(1) = \frac{\pi}{2}~.
\end{eqnarray}
From this we get:
\begin{eqnarray}
\sum^{\infty}_{k=1}\frac{1}{(2k+1)} 
\left(\begin{array}{c}
                               2k \\
                               k
\end{array}\right)\frac{1}{4^k}	
= \frac{\pi}{2} - 1~.
\end{eqnarray}
Analogously it follows: 
\begin{eqnarray}
\sum^{\infty}_{k=1}\frac{1}{(2k+1)^2}
\left(\begin{array}{c}
                               2k \\
                               k
\end{array}\right)\frac{1}{4^k}	
= \int^1_0 \frac{arcsin(x)}{x} dx -1 \nonumber 
\end{eqnarray}
\begin{eqnarray}
= \int^{\frac{\pi}{2}}_0 z cot(z) dz = \frac{\pi}{2}~ln(2) - 1~.
\end{eqnarray}
Therefore it follows:
\begin{eqnarray}
\sum^{\infty}_{k=1}\frac{1}{(2k+1)^3}
\left(\begin{array}{c}
                               2k \\
                               k
\end{array}\right)\frac{1}{4^k} = \int^1_0\frac{dy}{y} \left( \int^y_0 \frac{arcsin(x)}{x} \right) - 1~.
\end{eqnarray}
Substituting x=sin(z) and y = sin(p) gives us:
\begin{eqnarray}
\int^1_0\frac{dy}{y} \left( \int^y_0 \frac{arcsin(x)}{x} \right) =
	\int^{\frac{\pi}{2}}_0 \left( lnsin(p) \right)' dp \int^p_0 z \left( lnsin(z) \right)' dz~.
\end{eqnarray}
By partial integration it follows first:
\begin{eqnarray}
\int^1_0\frac{dy}{y} \left( \int^y_0 \frac{arcsin(x)}{x} \right) = 
- \int^{\frac{\pi}{2}}_0 p~lnsin(p)\left( lnsin(p) \right)' dp = \nonumber
\end{eqnarray}
\begin{eqnarray}
= -\frac{1}{2}\int^{\frac{\pi}{2}}_0 \left[ \left( lnsin(p) \right)^2 \right]' dp~.
\end{eqnarray}
A second partial integration results in:
\begin{eqnarray}
\int^1_0\frac{dy}{y} \left( \int^y_0 \frac{arcsin(x)}{x} \right) = \frac{1}{2} \int^{\frac{\pi}{2}}_0 \left( lnsin(p) \right)^2 dp~, 
\end{eqnarray}
thus we have:
\begin{eqnarray}
\sum^{\infty}_{k=1}\frac{1}{(2k+1)^3}
\left(\begin{array}{c}
                               2k \\
                               k
\end{array}\right)\frac{1}{4^k} 
= \frac{1}{2} \int^{\frac{\pi}{2}}_0 \left( lnsin(p) \right)^2 dp -1~.
\end{eqnarray}
Choosing n=3 in Eq.~(39) we can write:
\begin{eqnarray}
\sum^{\infty}_{k=1}\frac{1}{(2k+1)^4}
\left(\begin{array}{c}
                               2k \\
                               k
\end{array}\right)\frac{1}{4^k}	
= \int^1_0 \frac{dz}{z} \int^z_0 \frac{dy}{y} \int^y_0 \frac{dx}{x}  \frac{arcsin(x)}{x} dx -1 \nonumber 
\end{eqnarray}
Substituting x=sin(z), y = sin(t) and z = sin(p) and performing three partial integrations we get: 
\begin{eqnarray}
\sum^{\infty}_{k=1}\frac{1}{(2k+1)^4}
\left(\begin{array}{c}
                               2k \\
                               k
\end{array}\right)\frac{1}{4^k} 
= -\frac{1}{6} \int^{\frac{\pi}{2}}_0 \left( lnsin(p) \right)^3 dp -1~.
\end{eqnarray}
In order to calculate the n+1 order for the corresponding power series one has to perform n successive partial integrations
on the expression
\begin{eqnarray}
\sum^{\infty}_{k=1}\frac{1}{(2k+1)^{n+1}}
\left(\begin{array}{c}
                               2k \\
                               k
\end{array}\right)\frac{1}{4^k}	
= \int^1_0 \frac{dx_1}{x_1} \int^{x_1}_0~...~\int^{x_n}_0\frac{dx_n}{x_n} \frac{arcsin(x_n)}{x_n} dx_n -1 \nonumber 
\end{eqnarray}
Thus the theorem is proved.
\\
The integral on the right side is known for all n \cite{bow47,bor12}. For example, it follows from \cite{bor12}:
\begin{eqnarray}
Ls_4(\pi) = -\int^{\pi}_0 \left[ ln \left( 2sin \left(\frac{x}{2} \right) \right) \right]^3 dx~,
\end{eqnarray}
with
\begin{eqnarray}
Ls_4(\pi) = \frac{3}{2} \pi\zeta(3)~.
\end{eqnarray}
Substituting x = 2p it follows:
\begin{eqnarray}
\frac{3}{2} \pi\zeta(3) - \int^{\frac{\pi}{2}}_0 \left[ ln \left( 2sin \left( p \right) \right) \right]^3 dp = \nonumber
\end{eqnarray}
\begin{eqnarray}
- 2\int^{\frac{\pi}{2}}_0 \left[ ln(2) + lnsin \left( p \right) \right]^3 dp~. 
\end{eqnarray}
Therefore we get:
\begin{eqnarray}
	\int^{\frac{\pi}{2}}_0 \left[ lnsin \left( p \right) \right]^3 dp = -\frac{3}{4} \pi\zeta(3) - \frac{\pi^3}{8}ln(2) 
	-  \frac{\pi}{2}\left(ln(2)\right)^3~,
\end{eqnarray}
and the corresponding power series follows to:
\begin{eqnarray}
\sum^{\infty}_{k=1}\frac{1}{(2k+1)^4}
\left(\begin{array}{c}
                               2k \\
                               k
\end{array}\right)\frac{1}{4^k} = \frac{\pi}{8}\zeta(3) + \frac{\pi^3}{48}ln(2) + \frac{\pi}{12} \left( ln(2) \right)^3 - 1
\end{eqnarray}
As a further example we obtain for n=4:
\begin{eqnarray}
\sum^{\infty}_{k=1}\frac{1}{(2k+1)^5}
\left(\begin{array}{c}
                               2k \\
                               k
\end{array}\right)\frac{1}{4^k} = \frac{19\pi}{128}\zeta(4) + \frac{\pi}{8}ln(2)\zeta(3) + \frac{\pi^3}{96}\left( ln(2) \right)^2
+ \frac{\pi}{48}\left( ln(2) \right)^4 - 1~.
\end{eqnarray}

\section{A third kind of inverse power series with central Binomial coefficients and odd harmonic numbers}
The third family of series is defined as follows:

\begin{eqnarray}
v(n) = \sum^{\infty}_{k=1} \frac{h_k}{k^n} \left(\begin{array}{c}
                               2k \\
                               k
\end{array}\right)\frac{1}{4^k}~,
\end{eqnarray}
with $n \in \mathbb{N}$.

For n=1 the series is known from literature \cite{boy12}: 
\begin{eqnarray}
\sum^{\infty}_{k=1} \frac{h_k}{k} \left(\begin{array}{c}
                               2k \\
                               k
\end{array}\right)\frac{1}{4^k} = \frac{3}{2}\zeta(2)~.
\end{eqnarray}

\subsection{Lemma 3}
To be able to compute the corresponding series for higher n values we need a new identity again in form of a proper
valued integer series.

\begin{eqnarray}
g(k+1) = \sum^{\infty}_{i=1} \frac{1}{2k + 2i + 1}
\left(\begin{array}{c}
                               2i \\
                               i
\end{array}\right) \frac{1}{4^i} = \frac{\pi}{2}
\left(\begin{array}{c}
                               2k \\
                               k
\end{array}\right) \frac{1}{4^k} - \frac{1}{2k+1}~.
\end{eqnarray}

\subsection{Proof of Lemma 3}
For k=0 the series is known. Based on g(1) we can calculate for k=1 the corresponding value of g(2). The result is:
\begin{eqnarray}
\sum^{\infty}_{k=1} \frac{1}{2k+3} \left(\begin{array}{c}
                               2k \\
                               k
\end{array}\right)\frac{1}{4^k} =  \frac{\pi}{4} - \frac{1}{3}~.
\end{eqnarray}
Repeating the calculation k-times again an inhomogeneous difference equation of first order can be formulated. This
procedure has been introduced to prove lemma 1. It follows:
\begin{eqnarray}
g(k+1) - \frac{2k-1}{2k} g(k) = \frac{1}{2k(2k+1)}~.
\end{eqnarray}
The general solution of the inhomogeneous equation results to:
\begin{eqnarray}
g(k+1) = 
a \left(\begin{array}{c}
                               2k \\
                               k
\end{array}\right)\frac{1}{4^k} + 
\left(\begin{array}{c}
                               2k \\
                               k
\end{array}\right)\frac{1}{4^k}
\sum^k_{i=1} \frac{1}{2i(2i+1)}
\left(\begin{array}{c}
                               2i \\
                               i
\end{array}\right)^{-1}4^i~,
\end{eqnarray}
with an integer constant a. This constant follows from the boundary condition $g(1) = \frac{\pi}{2} -1$. The finite sum
appearing in the inhomogeneous part of the solution is known. It follows with \cite{spr06}:
\begin{eqnarray}
\sum^k_{i=1} \frac{1}{2i(2i+1)}
\left(\begin{array}{c}
                               2i \\
                               i
\end{array}\right)^{-1}4^i =
1 - \frac{1}{2k+1}
\left(\begin{array}{c}
                               2k \\
                               k
\end{array}\right)^{-1}4^k~.
\end{eqnarray}
From the boundary conditions we get:
\begin{eqnarray}
a = \frac{\pi}{2} -1~,	
\end{eqnarray}
and with this it follows:
\begin{eqnarray}
g(k+1) = \frac{\pi}{2}
\left(\begin{array}{c}
                               2k \\
                               k
\end{array}\right)\frac{1}{4^k} - \frac{1}{2k+1}~. 
\end{eqnarray}
Thus the lemma is proved.

Now we are able to calculate the corresponding series v(2). It is advantageous to start with g(k) instead with g(k+1). It follows:
\begin{eqnarray}
g(k) = \frac{\pi k}{2k-1}
\left(\begin{array}{c}
                               2k \\
                               k
\end{array}\right)\frac{1}{4^k} - \frac{1}{2k-1}~. 
\end{eqnarray}
First we calculate the expression
\begin{eqnarray}
\sum^{\infty}_{k=1} \frac{1}{(2k-1)^2} \sum^{\infty}_{i=1}\frac{1}{2i+2k-1}
\left(\begin{array}{c}
                               2i \\
                               i
\end{array}\right)\frac{1}{4^i} = \pi \sum^{\infty}_{k=1} \frac{k}{(2k-1)^3}
\left(\begin{array}{c}
                               2k \\
                               k
\end{array}\right)\frac{1}{4^k} - \sum^{\infty}_{k=1}\frac{1}{(2k-1)^3} \nonumber
\end{eqnarray}
\begin{eqnarray}
= \frac{\pi}{2} \sum^{\infty}_{k=0} \frac{k+1}{(2k+1)^3}
\left(\begin{array}{c}
                               2k+1 \\
                               k+1
\end{array}\right)\frac{1}{4^k} -\frac{7}{8}\zeta(3) \nonumber 
\end{eqnarray}
\begin{eqnarray}
= \frac{\pi}{2} \sum^{\infty}_{k=1} \frac{1}{(2k+1)^2}
\left(\begin{array}{c}
                               2k \\
                               k
\end{array}\right)\frac{1}{4^k} -\frac{7}{8}\zeta(3) + \frac{\pi}{2} \nonumber 
\end{eqnarray}
\begin{eqnarray}
= \frac{\pi}{2} \sum^{\infty}_{k=1} \frac{1}{(2k+1)^2}
\left(\begin{array}{c}
                               2k \\
                               k
\end{array}\right)\frac{1}{4^k} -\frac{7}{8}\zeta(3) + \frac{\pi}{2} \nonumber 
\end{eqnarray}
\begin{eqnarray}
= \frac{\pi}{2} \left( \frac{\pi}{2}ln(2) - 1 \right)  -\frac{7}{8}\zeta(3) + \frac{\pi}{2} \nonumber 
\end{eqnarray}
\begin{eqnarray}
= \frac{\pi^2}{4}ln(2) -\frac{7}{8}\zeta(3)~.
\end{eqnarray}
On the other hand it follows:
\begin{eqnarray}
\sum^{\infty}_{k=1} \frac{1}{(2k-1)^2} \sum^{\infty}_{i=1}\frac{1}{2i+2k-1}
\left(\begin{array}{c}
                               2i \\
                               i
\end{array}\right)\frac{1}{4^i} =  \sum^{\infty}_{i=1}
\left(\begin{array}{c}
                               2i \\
                               i
\end{array}\right)\frac{1}{4^i} \sum^{\infty}_{k=1} \frac{1}{(2k-1)^2(2i+2k-1)} \nonumber
\end{eqnarray}
\begin{eqnarray}
= \frac{1}{2}  \sum^{\infty}_{i=1}\frac{1}{i}
\left(\begin{array}{c}
                               2i \\
                               i
\end{array}\right)\frac{1}{4^i}\sum^{\infty}_{k=1}\frac{1}{(2k-1)^2} - \frac{1}{2}  \sum^{\infty}_{i=1}\frac{1}{i}
\left(\begin{array}{c}
                               2i \\
                               i
\end{array}\right)\frac{1}{4^i} \sum^{\infty}_{k=1} \frac{1}{(2k-1)(2i+2k-1)} \nonumber
\end{eqnarray}
\begin{eqnarray}
=\frac{3}{4}ln(2)\zeta(2) -\frac{1}{2}  \sum^{\infty}_{i=1}\frac{1}{i}
\left(\begin{array}{c}
                               2i \\
                               i
\end{array}\right)\frac{1}{4^i} \sum^{\infty}_{k=1} \frac{1}{(2k-1)(2i+2k-1)}~.
\end{eqnarray}
With
\begin{eqnarray}
\sum^{\infty}_{k=1} \frac{1}{(2k-1)(2i+2k-1)} = \frac{h_i}{2i}~,
\end{eqnarray}
where this identity follows immediately from the following partial fraction decomposition:
\begin{eqnarray}
\sum^{\infty}_{k=1} \frac{1}{(2k-1)(2i+2k-1)} = \frac{1}{2i} \sum^{\infty}_{k=1}
\left( \frac{1}{2k-1} - \frac{1}{2k+2i-1} \right)~.
\end{eqnarray}
We get:
\begin{eqnarray}
\sum^{\infty}_{k=1} \frac{1}{(2k-1)^2} \sum^{\infty}_{i=1}\frac{1}{2i+2k-1}
\left(\begin{array}{c}
                               2i \\
                               i
\end{array}\right)\frac{1}{4^i} = \frac{3}{4}ln(2)\zeta(2) - \frac{1}{4} \sum^{\infty}_{i=1}\frac{h_i}{i^2}
\left(\begin{array}{c}
                               2i \\
                               i
\end{array}\right)\frac{1}{4^i}~.
\end{eqnarray}
Comparing both sides and rearranging the different terms it follows:
\begin{eqnarray}
\sum^{\infty}_{k=1} \frac{h_k}{k^2} \left(\begin{array}{c}
                               2k \\
                               k
\end{array}\right)\frac{1}{4^k} = \frac{7}{2}\zeta(3) - 3ln(2)\zeta(2)~.
\end{eqnarray}
the corresponding series for higher n values can be obtained by calculating the following expression: 
\begin{eqnarray}
\sum^{\infty}_{k=1} \frac{1}{(2k-1)^n} \sum^{\infty}_{i=1}\frac{1}{2i+2k-1}
\left(\begin{array}{c}
                               2i \\
                               i
\end{array}\right)\frac{1}{4^i} 
\end{eqnarray}
As an example the series for n = 3, 4 and 5 have been calculated explicitly:
\begin{eqnarray}
\sum^{\infty}_{k=1} \frac{h_k}{k^3} \left(\begin{array}{c}
                               2k \\
                               k
\end{array}\right)\frac{1}{4^k} = \frac{15}{4}\zeta(4) - 7 ln(2)\zeta(3) + 3(ln(2))^2\zeta(2)
\end{eqnarray}

\begin{eqnarray}
\sum^{\infty}_{k=1}\frac{h_k}{k^4}
\left(\begin{array}{c}
                               2k \\
                               k
\end{array}\right)\frac{1}{4^k} = \frac{31}{2}\zeta(5) - \frac{15}{2}ln(2)\zeta(4) - \frac{13}{2}\zeta(2)\zeta(3) 
+ 7\left( ln(2) \right)^2\zeta(3) -2\left( ln(2) \right)^3\zeta(2)
\end{eqnarray}

\begin{eqnarray}
\sum^{\infty}_{k=1}\frac{h_k}{k^5}
\left(\begin{array}{c}
                               2k \\
                               k
\end{array}\right)\frac{1}{4^k} &=& \frac{399}{32}\zeta(6) - 31ln(2)\zeta(5) + \frac{15}{2}\left( ln(2) \right)^2\zeta(4)
+ 13ln(2)\zeta(2)\zeta(3) - 7\left( \zeta(3) \right)^2 - \frac{14}{3} \left( ln(2) \right)^3\zeta(3) \nonumber \\ 
&+& \left( ln(2) \right)^4\zeta(2)
\end{eqnarray}

\section{A fourth kind of inverse power series with central Binomial coefficients and odd harmonic numbers}
The fourth family of series defined as follows:

\begin{eqnarray}
z(n) = \sum^{\infty}_{k=1} \frac{h_k}{(2k-1)^n} \left(\begin{array}{c}
                               2k \\
                               k
\end{array}\right)\frac{1}{4^k}~,
\end{eqnarray}
with $n \in \mathbb{N}$. For n=1 we use Eq.~(10). The result is:
\begin{eqnarray}
\sum^{\infty}_{k=1} \frac{h_k}{2k-1}
\left(\begin{array}{c}
                               2k \\
                               k
\end{array}\right)\frac{1}{4^k} =
\frac{\pi}{2}~.
\end{eqnarray}

For higher n values the calculational scheme is as follows. Starting with Eq.~(68) we write:
\begin{eqnarray}
\sum^{\infty}_{i=1} \frac{1}{(2i-1)^2} 
\left(\begin{array}{c}
                               2i \\
                               i
\end{array}\right)\frac{1}{4^i}	
\left(\sum^{\infty}_{k=1} \frac{1}{(2k-1)(2k+2i-1)} \right) = \frac{1}{2} \sum^{\infty}_{i=1}
\frac{h_i}{i(2i-1)^2}
\left(\begin{array}{c}
                               2i \\
                               i
\end{array}\right)\frac{1}{4^i} \nonumber
\end{eqnarray}
\begin{eqnarray}
= \sum^{\infty}_{k=1} \frac{1}{2k-1} 
\left(\begin{array}{c}
                               2k \\
                               k
\end{array}\right)\frac{1}{4^k}	
\left(\sum^{\infty}_{i=1} \frac{1}{(2i-1)^2(2k+2i-1)}
\left(\begin{array}{c}
                               2i \\
                               i
\end{array}\right)\frac{1}{4^i} \right) \nonumber 
\end{eqnarray}
\begin{eqnarray}
= \frac{\pi}{2}ln(2) - \frac{1}{4}\sum^{\infty}_{k=1}\frac{1}{k^2(2k-1)} + \frac{1}{4} 
\sum^{\infty}_{k=1}\frac{1}{k^2(2k-1)^2}\left(\pi k 
\left(\begin{array}{c}
                               2k \\
                               k
\end{array}\right)\frac{1}{4^k} - 1 \right) \nonumber
\end{eqnarray}
\begin{eqnarray}
= \pi ln(2) - \pi + \frac{3}{4}\zeta(2)~.
\end{eqnarray}
The term on the right side of Eq.~(78) can be calculated by partial fraction decomposition, and from this we get:
\begin{eqnarray}
\sum^{\infty}_{k=1} \frac{h_k}{(2k-1)^2}
\left(\begin{array}{c}
                               2k \\
                               k
\end{array}\right)\frac{1}{4^k} =
\pi ln(2) - \frac{\pi}{2}~.
\end{eqnarray}
As an example for n=3 and n=4 it follows:
\begin{eqnarray}
	\sum^{\infty}_{k=1} \frac{h_k}{(2k-1)^3}
\left(\begin{array}{c}
                               2k \\
                               k
\end{array}\right)\frac{1}{4^k} =
\frac{3}{4} \pi \left( ln(2) \right)^2 - \pi ln(2) + \frac{\pi}{2}~,
\end{eqnarray}
\begin{eqnarray}
	\sum^{\infty}_{k=1} \frac{h_k}{(2k-1)^4}
\left(\begin{array}{c}
                               2k \\
                               k
\end{array}\right)\frac{1}{4^k} =
\frac{\pi}{16}\zeta(3) + \frac{\pi}{8}ln(2)\zeta(2) + \frac{\pi}{3} \left( ln(2) \right)^3
- \frac{3}{4}\pi \left( ln(2) \right)^2 + \pi ln(2) - \frac{\pi}{2}~.
\end{eqnarray}

\section{An application to nonlinear Euler sums}
In this section we present an application to nonlinear Euler sums which can be explicitly calculated in terms of the third kind
of power series which we have introduced in section V. These nonlinear Euler sums are defined as follows:
\begin{eqnarray}
w(n) = \sum^{\infty}_{k=1}\frac{h_k^2}{k^{2n}}~,
\end{eqnarray}
for $n \in \mathbb{N}$. A new identity is needed here:

\subsection{Lemma 4}
The following identity holds:
\begin{eqnarray}
\frac{h_k}{k} = 
\left(\begin{array}{c}
                               2k \\
                               k
\end{array}\right)\frac{1}{4^k}	
\sum^{\infty}_{k=1}\frac{h_i}{i+k}
\left(\begin{array}{c}
                               2i \\
                               i
\end{array}\right)\frac{1}{4^i}~.
\end{eqnarray}

\subsection{Proof of Lemma 4}
We start the proof with the calculation of the following expression:
\begin{eqnarray}
\sum^{\infty}_{i=1}\frac{h_i}{i+1}
\left(\begin{array}{c}
                               2i \\
                               i
\end{array}\right)\frac{1}{4^i} \nonumber
\end{eqnarray}
\begin{eqnarray}
= \sum^{\infty}_{i=1}\frac{h_{i+1}}{i+1}
\left(\begin{array}{c}
                               2i \\
                               i
\end{array}\right)\frac{1}{4^i} - \sum^{\infty}_{i=1} \frac{1}{(i+1)(2i+1)}
\left(\begin{array}{c}
                               2i \\
                               i
\end{array}\right)\frac{1}{4^i} \nonumber
\end{eqnarray}
\begin{eqnarray}
= 4\sum^{\infty}_{i=2}\frac{h_{i}}{i}
\left(\begin{array}{c}
                               2i-2 \\
                               i-1
\end{array}\right)\frac{1}{4^i} - \sum^{\infty}_{i=1} \frac{1}{(i+1)(2i+1)}
\left(\begin{array}{c}
                               2i \\
                               i
\end{array}\right)\frac{1}{4^i} \nonumber
\end{eqnarray}
\begin{eqnarray}
= 2\sum^{\infty}_{i=1}\frac{h_{i}}{2i-1}
\left(\begin{array}{c}
                               2i \\
                               i
\end{array}\right)\frac{1}{4^i} - \sum^{\infty}_{i=1} \frac{1}{(i+1)(2i+1)}
\left(\begin{array}{c}
                               2i \\
                               i
\end{array}\right)\frac{1}{4^i} - 1.
\end{eqnarray}
Using Eq.(10) with $h(i)= h_i$ it follows:
\begin{eqnarray}
\sum^{\infty}_{i=1}\frac{h_{i}}{2i-1}
\left(\begin{array}{c}
                               2i \\
                               i
\end{array}\right)\frac{1}{4^i} = 1 + \sum^{\infty}_{i=1}\frac{1}{2i+1}
\left(\begin{array}{c}
                               2i \\
                               i
\end{array}\right)\frac{1}{4^i}~,
\end{eqnarray}
and therefore we get: 
\begin{eqnarray}
\sum^{\infty}_{i=1}\frac{h_i}{i+1}
\left(\begin{array}{c}
                               2i \\
                               i
\end{array}\right)\frac{1}{4^i} = 2
\end{eqnarray}
Now we perform the following calculation: 
\begin{eqnarray}
\sum^{\infty}_{i=1}\frac{h_i}{i+1}
\left(\begin{array}{c}
                               2i \\
                               i
\end{array}\right)\frac{1}{4^i} \nonumber
\end{eqnarray}
\begin{eqnarray}
=\sum^{\infty}_{i=0}\frac{h_{i+1}}{i+2}
\left(\begin{array}{c}
                               2i+2 \\
                               i+1
\end{array}\right)\frac{1}{4^{i+1}} \nonumber 
\end{eqnarray}
\begin{eqnarray}
=\frac{1}{2}\sum^{\infty}_{i=0}\frac{(2i+1)h_{i+1}}{(i+1)(i+2)}
\left(\begin{array}{c}
                               2i \\
                               i
\end{array}\right)\frac{1}{4^i} \nonumber 
\end{eqnarray}
\begin{eqnarray}
=\sum^{\infty}_{i=1}\frac{h_i}{i+2}
\left(\begin{array}{c}
                               2i \\
                               i
\end{array}\right)\frac{1}{4^i} -
\frac{1}{2}\sum^{\infty}_{i=1}\frac{h_i}{(i+1)(i+2)}
\left(\begin{array}{c}
                               2i \\
                               i
\end{array}\right)\frac{1}{4^i} +
\frac{1}{2}\sum^{\infty}_{i=0}\frac{1}{(i+1)(i+2)}
\left(\begin{array}{c}
                               2i \\
                               i
\end{array}\right)\frac{1}{4^i}~.
\end{eqnarray}
Rearranging the different terms we find:
\begin{eqnarray}
\sum^{\infty}_{i=1}\frac{h_i}{i+2}
\left(\begin{array}{c}
                               2i \\
                               i
\end{array}\right)\frac{1}{4^i} = 
\frac{2}{3}\sum^{\infty}_{i=1}\frac{h_i}{i+1}
\left(\begin{array}{c}
                               2i \\
                               i
\end{array}\right)\frac{1}{4^i} - \frac{1}{3} 
\sum^{\infty}_{i=0}\frac{1}{(i+1)(i+2)}
\left(\begin{array}{c}
                               2i \\
                               i
\end{array}\right)\frac{1}{4^i} + \frac{2}{3}~.
\end{eqnarray}
The analogous calculation for k=3 results in:
\begin{eqnarray}
\sum^{\infty}_{i=1}\frac{h_i}{i+3}
\left(\begin{array}{c}
                               2i \\
                               i
\end{array}\right)\frac{1}{4^i} = 
\frac{4}{5}\sum^{\infty}_{i=1}\frac{h_i}{i+2}
\left(\begin{array}{c}
                               2i \\
                               i
\end{array}\right)\frac{1}{4^i} - \frac{2}{5} 
\sum^{\infty}_{i=0}\frac{1}{(i+1)(i+3)}
\left(\begin{array}{c}
                               2i \\
                               i
\end{array}\right)\frac{1}{4^i} + \frac{2}{5}~.
\end{eqnarray}
Repeating this procedure k times again a inhomogeneous difference equation can be formulated:
\begin{eqnarray}
o(k) - \frac{2k-2}{2k-1} o(k-1) = p(k)~, 	
\end{eqnarray}
where the solution of the homogeneous solutions results to: 
\begin{eqnarray}
o_h(k) = \frac{1}{k}
\left(\begin{array}{c}
                               2k \\
                               k
\end{array}\right)^{-1}4^k~.
\end{eqnarray}
The inhomogeneity results by partial fraction decomposition and use of lemma 1 to:
\begin{eqnarray}
p(k) = \frac{2}{2k-1} - \frac{1}{2k-1} + \frac{1}{k(2k-1)}
\left(\begin{array}{c}
                               2k \\
                               k
\end{array}\right)^{-1}4^k - \frac{1}{k(2k-1)} -  \frac{k-1}{k(2k-1)} \nonumber
\end{eqnarray}
\begin{eqnarray}
= \frac{1}{k(2k-1)}
\left(\begin{array}{c}
                               2k \\
                               k
\end{array}\right)^{-1}4^k~.  
\end{eqnarray}
Therefore the solution of the corresponding difference equation is:
\begin{eqnarray}
	p(k) = \frac{1}{k}
\left(\begin{array}{c}
                               2k \\
                               k
\end{array}\right)^{-1}4^k \sum^k_{i=1} i
\left(\begin{array}{c}
                               2i \\
                               i
\end{array}\right)\frac{1}{4^i} \frac{1}{i(2i-1)}
\left(\begin{array}{c}
                               2i \\
                               i
\end{array}\right)^{-1}4^i \nonumber
\end{eqnarray}
\begin{eqnarray}
= \frac{h_k}{k}
\left(\begin{array}{c}
                               2k \\
                               k
\end{array}\right)^{-1}4^k~. 
\end{eqnarray}
Thus lemma 4 is proved.
\\
Now we are able to calculate the nonlinear Euler sum for n=1. We start with the expression:  
\begin{eqnarray}
\sum^{\infty}_{k=1} \frac{h^2_k}{k^2} = \sum^{\infty}_{k=1} \frac{h_k}{k} 
\left(\begin{array}{c}
                               2k \\
                               k
\end{array}\right)\frac{1}{4^k}
\sum^{\infty}_{i=1}\frac{h_i}{i+k}
\left(\begin{array}{c}
                               2i \\
                               i
\end{array}\right)\frac{1}{4^i}   \nonumber
\end{eqnarray}
\begin{eqnarray}
= \sum^{\infty}_{i=1} h_i
\left(\begin{array}{c}
                               2i \\
                               i
\end{array}\right)\frac{1}{4^i}
\sum^{\infty}_{k=1}\frac{h_k}{k(i+k)}
\left(\begin{array}{c}
                               2k \\
                               k
\end{array}\right)\frac{1}{4^k}   \nonumber
\end{eqnarray}
\begin{eqnarray}
= \sum^{\infty}_{i=1} \frac{h_i}{i}
\left(\begin{array}{c}
                               2i \\
                               i
\end{array}\right)\frac{1}{4^i}
\sum^{\infty}_{k=1}\frac{h_k}{k}
\left(\begin{array}{c}
                               2k \\
                               k
\end{array}\right)\frac{1}{4^k}  -
\sum^{\infty}_{i=1} \frac{h_i}{i}
\left(\begin{array}{c}
                               2i \\
                               i
\end{array}\right)\frac{1}{4^i}
\sum^{\infty}_{k=1}\frac{h_k}{i+k}
\left(\begin{array}{c}
                               2k \\
                               k
\end{array}\right)\frac{1}{4^k}~. 
\end{eqnarray}
Thus it follows:
\begin{eqnarray}
\sum^{\infty}_{k=1} \frac{h^2_k}{k^2} = \left(
\sum^{\infty}_{k=1} \frac{h_k}{k}
\left(\begin{array}{c}
                               2k \\
                               k
\end{array}\right)\frac{1}{4^k} \right)^2 = \frac{9}{8} \left( \zeta(2) \right)^2 \nonumber
\end{eqnarray}
\begin{eqnarray}
= \frac{45}{16} \zeta(4)~.
\end{eqnarray}
For higher n values we calculate the following expression: 
\begin{eqnarray}
\sum^{\infty}_{k=1} \frac{h_k}{k^{(2n-1)}} 
\left(\begin{array}{c}
                               2k \\
                               k
\end{array}\right)\frac{1}{4^k}
\sum^{\infty}_{i=1}\frac{h_i}{i+k}
\left(\begin{array}{c}
                               2i \\
                               i
\end{array}\right)\frac{1}{4^i} 
\end{eqnarray}
For example, for n=2 and 3 it follows:
\begin{eqnarray}
\sum^{\infty}_{k=1}\frac{h_k^2}{k^4} = \frac{315}{32}\zeta(6) - \frac{49}{8}\zeta(3)^2~,
\end{eqnarray}
\begin{eqnarray}
\sum^{\infty}_{k=1}\frac{h_k^2}{k^6} = \frac{315}{8}\zeta(8) - \frac{217}{4}\zeta(3)\zeta(5) + \frac{49}{4}\zeta(2)\zeta(3)^2~.
\end{eqnarray}
Thus for all even powers of k in the denominator the concerning Euler sums can be calculated explicitly and purely in terms
of zeta functions.

\section{outlook}
For odd powers the situation seems to be similar as we have: 

\subsection{Lemma 5}
\begin{eqnarray}
\sum^{\infty}_{k=1}\frac{h_k^2}{k^3} = \frac{7}{4}\zeta(2)\zeta(3) - \frac{31}{16}\zeta(5)~.
\end{eqnarray}
This is different from the case where $h_k$ appears in the nominator instead of $h_k^2$. Here it holds \cite{sit85}: 
\begin{eqnarray}
\sum^{\infty}_{k=1}\frac{h_k}{k^3} = - \frac{53}{8}\zeta(4) + 7 ln(2) \zeta(3) - 2\left( ln(2)\right)^2 \zeta(2)
+ \frac{1}{3}\left( ln(2)\right)^4 + 8Li_4\left(\frac{1}{2}\right)~,
\end{eqnarray}
where the Polylogarithm function appears.  

\subsection{Proof of Lemma 5}
Eq.~(99) is obtained with the help of the following two equalities which are known from literature \cite{val19}:
\begin{eqnarray}
\sum^{\infty}_{k=1}\frac{H_k H_{2k}}{(2k)^3}&=&\frac{307}{128}\zeta(5)-\frac{1}{16}\zeta(2)\zeta(3)+\frac{1}{3}\left( ln(2)\right)^3 \zeta(2)
-\frac{7}{8}\left( ln(2) \right)^2\zeta(3)- \frac{1}{15}\left( ln(2)\right)^5 \nonumber \\
&-& 2ln(2)Li_4\left(\frac{1}{2}\right)- 2Li_5\left(\frac{1}{2}\right)~,
\end{eqnarray}
where $Li_n$ denotes the Polylogarithm function. The second equality holds:
\begin{eqnarray}
\sum^{\infty}_{k=1}(-1)^{k-1}\frac{H_k^2}{k^3}&=&\frac{2}{15}\left( ln(2)\right)^5 - \frac{11}{8}\zeta(2)\zeta(3)-\frac{19}{32}\zeta(5)
+\frac{7}{4}\left( ln(2)\right)^2 \zeta(3) -\frac{2}{3}\left( ln(2)\right)^3 \zeta(2) \nonumber \\
&+& 4ln(2)Li_4\left(\frac{1}{2}\right)+ 4Li_5\left(\frac{1}{2}\right)~.
\end{eqnarray}
From Eq.~(101) it follows by decomposing the corresponding sum in even and odd contributions:
\begin{eqnarray}
\sum^{\infty}_{k=1}\frac{H_k h_k}{k^3} = \frac{279}{16}\zeta(5) - 7\left( ln(2)\right)^2 \zeta(3)
+\frac{8}{3}\left( ln(2)\right)^3 \zeta(2) - \frac{8}{15}\left( ln(2)\right)^5
- 16ln(2)Li_4\left(\frac{1}{2}\right)- 16Li_5\left(\frac{1}{2}\right)~.
\end{eqnarray}
Now we write
\begin{eqnarray}
\sum^{\infty}_{k=1}\frac{H_k^2}{k^3} - \sum^{\infty}_{k=1}(-1)^{k-1}\frac{H_k^2}{k^3} = 
\frac{1}{4} \sum^{\infty}_{k=1}\frac{{H_{2k}^2}}{k^3}~. 
\end{eqnarray}
From this equality we obtain by rearranging the different terms:
\begin{eqnarray}
\sum^{\infty}_{k=1}\frac{h_k^2}{k^3} = \frac{15}{4}\sum^{\infty}_{k=1}\frac{H_k^2}{k^3} - 4\sum^{\infty}_{k=1}(-1)^{k-1}\frac{H_k^2}{k^3}
-\sum^{\infty}_{k=1}\frac{H_k h_k}{k^3} = \frac{7}{4}\zeta(2)\zeta(3) - \frac{31}{16}\zeta(5)~.
\end{eqnarray}
Thus Lemma 5 is proved, and as a consequence we suppose that the corresponding Euler sums with higher odd powers of k in the denominator
can also be expressed purely in terms of zeta function.

\section{Summary}
We have introduced a special summation method that allows, based on proper two-valued integer functions, to calculate explicitly
in terms of natural logarithms and zeta functions a variety of power series that involves central binomial numbers and in
addition odd harmonic numbers. With the help of these series we have shown that a special type of nonlinear Euler sums, where
explicitly inverse power sums of even degree are combined with odd harmonic numbers of second degree, can be calculated
explicitly and purely in terms of zeta functions.

\end{document}